\documentclass[12pt]{article}
\title{Efficiency for multitime vector variational problems
    on Riemannian manifolds involving geodesic quasiinvex functionals}
\author{\c{S}tefan Mititelu, M\u ad\u alina Constantinescu, Constantin Udri\c{s}te}
\date{}
\usepackage{graphicx}
\begin{document}
\maketitle

\newtheorem{Th}{Theorem}
\newtheorem{Co}{Corollary}
\newtheorem{Prop}{Proposition}
\newtheorem{Lem}{Lemma}
\newtheorem{Def}{Definition}

\begin{abstract}
We study the connection between a multitime scalar
variational problem (SVP), a multitime vector variational problem (VVP) and
a multitime vector fractional variational problem (VFP). For (SVP),
we establish necessary optimality conditions. For both vector variational
problems, we define the notions of Pareto efficient solution
and of normal efficient solution and we establish necessary efficiency
conditions for (VVP) and (VFP) using both notions. The main
purpose of the paper is to establish sufficient efficiency
conditions for the vector problems (VVP) and (VFP). Moreover,
we obtain sufficient optimality conditions for (SVP).
The sufficient conditions are based on our original notion of
$(\rho ,b)$-geodesic quasiinvexity.
\end{abstract}



\newcommand{\omm}{\Omega}
\newcommand{\om}{\omega}
\newcommand{\ty}{\infty}
\newcommand{\di}{\displaystyle}
\newcommand{\va}{\varphi}
\newcommand{\si}{\sigma}
\newcommand{\ga}{\gamma}
\newcommand{\gaa}{\Gamma}
\newcommand{\na}{\nabla}
\newcommand{\te}{\theta}
\newcommand{\ld}{\ldots}
\newcommand{\ov}{\over}
\newcommand{\ri}{\Rightarrow}
\newcommand{\rii}{\Leftrightarrow}
\newcommand{\noa}{\noalign{\medskip}}
\newcommand{\mm}{\medskip}
\newcommand{\la}{\lambda}
\newcommand{\su}{\subset}
\newcommand{\qu}{\quad}
\newcommand{\fo}{\forall}
\newcommand{\al}{\alpha}
\newcommand{\be}{\beta}
\newcommand{\ep}{\varepsilon}
\newcommand{\pa}{\partial}
\newcommand{\ti}{\times}
\newcommand{\de}{\delta}
\newcommand{\dd}{\hbox{d}}
\newcommand{\pp}{\mapsto}
\newcommand{\tg}{\hbox{tg}}

\section{Introduction and preliminaries}

Beginning with Valentine [16] in 1937, during the years, the variational problem with constraints
experienced different stages of development. In 2008, Mititelu [6] studied a single-time vector
(or multiobjective) fractional variational problem.
Mititelu [6] and Mititelu {\&} Stancu-Minasian [9] established for this problem necessary
efficiency (Pareto minimum) conditions. Using generalized quasiinvex functions,
they developed a duality theory including weak, direct and converse duality theorems.

In 2007 Udri\c{s}te and \c{T}evy [15] gave new results for a multitime variational problem of vector variable.
In 2009 Pitea, Udri\c{s}te and Mititelu [11], [12] considered the multitime vector variant of the problem
(MSP) into geometrical language, using curvilinear integrals, and establishing necessary efficiency conditions
and developed a duality theory for this problem. Recently Mititelu and Postolache [8] studied the
same subject for multitime vector fractional and nonfractional variational problems on
Riemannian manifolds, but using multiple integrals.

The purpose of this work is to deduce necessary and sufficient optimality conditions for
the multitime scalar problem (SVP) (sections 2, 4) and of Pareto efficiency for the multitime
vector variational problems (VVP) and (VFP) (section 3, 4), in a geometrical framework [11], [12].

Let $(T,h) $and $(M,g)$ be two Riemannian manifolds of dimensions $m$ and $n$.
In addition, $M$ is a complete manifold. Denote $t = (t^1,...,t^m) = (t^\upsilon)$
the points of a measurable set $\Omega$ in $T$ and $x = (x^1,...,x^n) =
(x^i)$ the points of $M$. Consider the first order jet bundle $J^1(T,M) = \Omega \times
R^n \times R^{nm}$ and the functions
$$
x:\Omega\subset T \to M,\,\, X:J^1(T,M) \to \mbox{R},
$$
$$
f = (f_r ): J^1(T,M) \to R^p, \; k = (k_r ):J^1(T,M) \to R^p,
$$
$$
g = (g_\alpha ):J^1(T,M) \to \mbox{R}^m, \;
h = (h_s ):J^1(T,M) \to R^q,
$$
where $m, p, q \in N^\ast$, $r = \overline {1,p}$, $\alpha = \overline {1,m}$ and $s = \overline {1,q} $, all of $C^2$-class.

The argument of each function $X,f,k,g,h$ is $j^1x=(t,x,x_\upsilon)$, the first prolongation jet of $x$. For functionals, based on Lagrangians $X,f,k,g,h$,
 we use the pullback $j^1_tx=(j^1x)(t))$, where
$t \in \Omega$ and $x(t) = (x^i(t))$, and $(\partial x / \partial t^\upsilon (t)) = (x_\upsilon(t)) = (x^i_\upsilon(t))$.

The Euler-Ostrogradsky PDEs produced by the Lagrangian $X$ are
$$
\frac{\partial X}{\partial x^k} - \frac{\partial }{\partial
t^\upsilon }\left( {\frac{\partial X}{\partial x^k_\upsilon}} \right) = 0, \, k=\overline{1,n};\,\upsilon = \overline {1,m}.\leqno (1.1)
$$

To develop our theory, we shall use a normed vector space of functions $(\mathcal{F}(\Omega,M), ||\cdot||)$, where
\[
{\mathcal{F}(\Omega,M)} = \{x:\Omega \to M\;\left| {\;x}
\right.\;\;\mbox{is}\;\mbox{piecewice}\;C^1\},
\]
and
\[
\left\| x \right\| = \left\| x \right\|_\infty + \sum^m_{\ga =1} \sum^n_{k=1} ||x^k_\ga||_\infty.
\]
The induced distance is $d(x^0(\cdot), x(\cdot))=||x(\cdot)- x^0(\cdot)||$, $x^0(\cdot)$, $x(\cdot)\in \mathcal{F}(\Omega,M)$. In this sense,
$(\mathcal{F}(\Omega,M), d)$ is a metric space.

The following partial ordering is used for two $n$-tuples $v = (v_1 ,\ldots ,v_n )$ and $w = (w_1 , \ldots, w_n )$:
\[
\begin{array}{ll}
v = w\; \Leftrightarrow \;v_i = w_i, \; i = \overline {1,n} &
v < w\; \Leftrightarrow \;v_i < w_i , \;i = \overline {1,n} \\ \noa
v \preceq w \Leftrightarrow v_i \leq w_i, \; i = \overline {1,n}; & v \le w \Leftrightarrow v \preceq w \; \hbox{and} \; v \ne w. \end{array}
$$
Similar partial relations are used also for $m$-tuples.

Let $dv = \sqrt{\det h}\,\,dt^1\wedge dt^2 \ldots \wedge dt^m$ be the volume element on $\Omega$. We use the functionals
$$F_r(x(\cdot))=\int_a^b f_r (j^1_tx)dv, \,\,K_r(x(\cdot))=\int_a^b k_r (j^1_tx)dv,\,r=\overline{1,p}$$
and a vector fractional functional
$$J(x(\cdot))=\left( \frac{F_1}{K_1 }(x(\cdot)) , \ldots , \frac{F_p }{K_p }(x(\cdot)) \right).$$
The general problem of study is the multitime vector fractional problem
$$
\hbox{(VFP)}\; \left\{ {\begin{array}{l}
\mbox{Maximize}\;\mbox{Pareto}\; J(x(\cdot))\\ \noa
 \mbox{subject}\;\mbox{to}\quad \quad g(j^1_tx)\;\preceq\;0,\; h(j^1_tx)\; = \;0,\; \\
 \quad \quad \quad \quad \quad \quad x(t)\left| {_{\partial \Omega } } \right. = u(t)\;(\mbox{given}),\;\forall t \in \Omega.
 \end{array}} \right.
$$
\noindent

This problem include the multitime vector variational problem
$$
\hbox{(VVP)}\; \left\{ {\begin{array}{l}
\mbox{Minimize}\;\mbox{Pareto} \; \left( F_1(\cdot), \ld ,\di F_p(x(\cdot) \right) \\ \noa
 \mbox{subject}\;\;\mbox{to}\quad g(j^1_tx)\;\preceq\;0, \; h(j^1_tx)\; = 0,\; \\
 \quad \quad \quad \quad \quad \;x(t)\left| {_{\partial \Omega } =
u(t),\;\forall t \in \Omega \;.\;} \right. \\
 \end{array}} \right.
$$
and the following multitime scalar variational problem
$$
\left\{ {\begin{array}{l}
\mbox{Minimize}\; E(x(\cdot)) = \int_\Omega {X(j^1_tx)dv} \\ \noa  \mbox{subject}\;\;\mbox{to}\;\;g(j^1_tx)\;\preceq\;0; \; h(j^1_tx)\; = \;0,\;t \in \Omega \\
 \quad \quad \quad \quad \;x(t)\left| {_{\partial \Omega } = u(t).} \right.  \end{array}} \right. \leqno \hbox{(SVP)}
$$
The three  variational problems have the same domain
{\small
$$
\mathcal{D} = \left\{ {x \in \mathcal{F}(\Omega,M) } \right.\left| {g(j^1_tx) }
\right.\preceq\,0,\;h(j^1_tx) = 0,\;x(t)\left| {_{\partial \Omega } } \right. = u(t)\}.
$$}

\begin{Def} \cite{[1]} Let $(M,g)$ be a complete Riemannian manifold.
Let $\eta:M\times M \to TM, \,\eta(u,x)\in T_uM,\,\, u,x\in M$ be a vector function and
$S\subset M$ a nonempty set.

(ii) The set $S$ is called $\eta$-geodesic invex if, for every $u,x\in S$,
there exists exactly one geodesic $\gamma_{u,x}:[0,1]\to M$ such that
$$\gamma_{u,x}(0) =u, \dot \gamma_{u,x}(0)=\eta(u,x),\,\,\gamma_{u,x}(\tau)\in S, \,\,\forall \tau \in [0,1].$$

(ii) Let $S\subset M$ be an open $\eta$-geodesic invex set and $f:S\to R$ be a $C^1$ function.
The function $f$ is called $\eta$-geodesic invex on $S$ if
$$f(x)-f(u)\geq df_u(\eta(u,x)),\,\,\forall u,x \in S.$$
\end{Def}

\begin{Def} \cite{[16]} Let $x^0(\cdot), x(\cdot) \in \mathcal{F}(\Omega,M)$. A function
$\varphi(t,\tau)$, $t\in \Omega$, $\tau \in [0,1]$ is called geodesic deformation of the pair of
functions $(x^0(\cdot), x(\cdot))$, if it satisfies the properties: (1) the function $\tau \to \varphi(t,\tau)$ is a geodesic;
(2) $\varphi(t,0)=x^0(t), \varphi(t,1)=x(t)$.
\end{Def}

\begin{Def} The set ${\cal S}= {\cal F}(\Omega,S)\subset \mathcal{F}(\Omega,M)$ is called $\eta$-geodesic invex if, for every $x^0(\cdot), x(\cdot) \in {\cal S}$,
there exists exactly one geodesic deformation $\varphi(t,\tau)$, $t\in \Omega$, $\tau \in [0,1]$ such that the vector function
$$\eta(x^0(t), x(t))=\frac{\partial \varphi}{\partial \tau}(t,\tau)|_{\tau =0} \in T_{x^0(t)}M \equiv \eta (t) = (\eta^1 (t), \ld, \eta^n (t))$$
is of class $C^1$ and satisfies $\eta (t)\left| {_{\partial \Omega } } \right. = 0$.
\end{Def}

For our sufficient conditions of efficiency and optimality, we shall introduce the notion of $(\rho ,b)$-geodesic quasiinvex
functionals. We fix a number $\rho \in R$, a functional $b:\mathcal{F}(\Omega,M) \times \mathcal{F}(\Omega,M)  \to \;[0,\;\infty )$ and the
distance function $d(x(\cdot),y(\cdot))$ on $\mathcal{F}(\Omega,M)$. We consider the functional
\[
E:\mathcal{F}(\Omega,M) \to R, \,\,\,E(x(\cdot)) = \int_\Omega {X(j^1_tx)\;\mbox{dv}}.
\]

\begin{Def} Let $(M,g)$ be a complete Riemannian manifold.
Let ${\cal S}$ be an open $\eta$-geodesic invex subset of $\mathcal{F}(\Omega,M)$.

(i) The functional $E$ is called (strictly) $(\rho ,b)$-geodesic quasiinvex
at $x^0(\cdot)\in {\cal S}$, with respect to $\eta (t)$, if
$E(x(\cdot)) \leq\;E(x^0(\cdot))$ implies
\[
b(x,x^0)\int_\Omega {\left( {\eta^i \,\frac{\partial X}{\partial
x^i}(j^1_tx^0) + \frac{\partial \eta^i }{\partial
t^\upsilon }\frac{\partial X}{\partial x_\upsilon ^i }(j^1_tx^0)} \right)\;\mbox{d}v\;}
\;\;( < )\;\leq \,\,\,
 - \rho b(x,x^0)d^2(x,x^0),
\]
for any $x(\cdot)\in {\cal S}$.

(ii) The functional $E$ is called monotonic $(\rho ,b)$-geodesic quasiinvex
at $x^0(t)\in S$,  with respect to $\eta (t)$, if
$E(x(\cdot)) \leq\;E(x^0(\cdot))$ implies
\[
b(x,x^0)\int_\Omega {\left( {\eta^i\, \frac{\partial X}{\partial
x^i}(j^1_tx^0) + \frac{\partial \eta^i }{\partial t^\upsilon }\frac{\partial X}{\partial x_\upsilon ^i }(j^1_tx^0)} \right)\;\mbox{d}v\;  }
= - \rho b(x,x^0)d^2(x,x^0),
\]
for any $x(\cdot)\in {\cal S}$.
\end{Def}

{\bf Example} Let us fix the domain
$$E =\{ x: \Omega = [0,1]^m\subset R^m \to R_+\,\,\Big|\,\,x(\cdot)\,\, \hbox{continuous}\}$$
and the "negative" Boltzmann-Shannon functional
$$J:E\to R,\,\, J(x(\cdot))= \int_\Omega x(t)\ln x(t)dv.$$
This functional is geodesic quasiinvex with respect to
$$\eta(t)=\left\{\begin{array}{ccc} - (\ln x(t)+1)\,d^2(x,x^0)& if & t\in \hbox{int}\,\, \Omega\\ \
0&if& t\in \partial \Omega.\end{array}\right.$$

\section{Necessary optimality conditions for \\scalar problem (SVP)}

In all the paper, we simplify supposing $T=R^m$ and hence $\det h = 1$.

We start with variational problem (SVP), recalling a well-known result.

\begin{Th}
(Necessary optimality to (SVP) (Mititelu, Postolache [8, Theorem 2.1]). If $x^0(\cdot) \in \mathcal{D}$ is an optimal solution to problem (SVP), then there exist a scalar $\tau \in R$ and the piecewise smooth multipliers $\lambda (t) = (\lambda^\alpha (t)) \in R^m$ and $\mu (t) = (\mu^s (t)) \in R^q$ that satisfy the following conditions:
$$
\left\{ {\begin{array}{l}
\tau \di\frac{\partial X}{\partial x^0} + \lambda^\alpha (t) \di\frac{\partial g_\alpha }{\partial x^0} + \mu^s (t) \di\frac{\partial h_s }{\partial x^0} - \\ \noa
- \di\frac{\partial}{\partial t^\upsilon } \left( {\tau \di\frac{\partial X}{\partial x_\upsilon^0 } + \lambda^\alpha (t) \di\frac{\partial g_\alpha}{\partial x_\upsilon^0 } + \mu^s (t) \frac{\partial h_s }{\partial x_\upsilon^0}} \right) = 0 \\ \noa
\lambda^\alpha (t) g_\alpha (j^1_tx^0) = 0, \; \hbox{for each $\alpha = \overline{1,m}$ (no summation)} \\ \noa
\tau \;\geq 0, \quad (\lambda ^\alpha (t))\;\succeq\;0,\;t \in \Omega, \end{array}} \right. \leqno \hbox{(SFJ)}
$$
where $$x^0 = (x^k)^0, \di\frac{\partial X}{\partial x^0}:  = \di\frac{\partial X}{\partial x}(j^1_tx^0), \di\frac{\partial X}{\partial x_\upsilon^0 } : = \di\frac{\partial X}{\partial x_\upsilon }(j^1_tx^0)$$ etc.
\end{Th}

\begin{Def}
A point $x^0(\cdot) \in \mathcal{D}$ is called normal optimal solution to (SVP) if $\tau > 0$.
\end{Def}

\section{Necessary efficiency conditions for \\multitime vector variational problems \\(VVP) and (VFP)}

\subsection{Efficiency for multitime vector variational \\problems (VVP)}

We consider the vector functional
$$
F(x(\cdot)) = \left( {F_1 (x(\cdot)), \ld, F_p (x(\cdot))} \right)
$$
and the multitime vector variational problem
$$
\hbox{(VVP)}\; \left\{ {\begin{array}{l}
\mbox{Minimize Pareto}\qu F(x(\cdot)) \\ \noa
\mbox{subject to} \quad g(j^1_tx)\;\preceq\;0, \; h(j^1_tx)\; = 0,\; \\
 \quad \quad \quad \quad \quad \;x(t)\left| {_{\partial \Omega } =
u(t),\;\forall t \in \Omega \;.\;} \right. \\
 \end{array}} \right.
$$
The domain of (VVP) is just $\mathcal{D}$.

In this section we establish necessary Pareto efficiency conditions for the program (VVP).

\begin{Def}
 A function $x^0(\cdot) \in \mathcal{D}$ is called an efficiency solution
(Pareto minimum) to (VVP) if there exists no $x(\cdot) \in \mathcal{D}$ such that $F(x(\cdot)) \preceq F(x^0(\cdot))$.
\end{Def}

\begin{Th}
(Necessary efficiency to (VVP) (Mititelu, Postolache [8, Theorem 3.1]). Consider the vector multitime variational problem (VVP) in the framework presented in Section 1.1 and let $x^0(\cdot) \in \mathcal{D}$ be an efficiency solution to (VVP). Then there are the vector Lagrange multipliers $\tau = (\tau^r ) \in {\bf R}^p$ and $\lambda (t) = (\lambda^\alpha (t)) \in R^m$ and $\mu (t) = (\mu^s (t)) \in R^q$, all functions being piecewise smooth, which satisfy the conditions
$$
\left\{ \begin{array}{l}
\tau^r \di\frac{\partial f_r}{\partial x^0} + \lambda^\alpha
(t) \di\frac{\partial g_\alpha }{\partial x^0} + \mu^s (t) \di\frac{\partial h_s }{\partial x^0} - \\ \noa
- \di\frac{\partial }{\partial t^\upsilon }\left( {\tau^r \di\frac{\partial f_r }{\partial x_\upsilon^0 } + \lambda^\alpha (t) \di\frac{\partial g_\alpha }{\partial x_\upsilon^0 } + \mu^s
(t) \di\frac{\partial h_s }{\partial x_\upsilon^0 }} \right) = 0 \\ \noa
\lambda^\alpha (t) g_\alpha (j^1_tx^0) = 0, \qu \hbox{for each} \qu \alpha = \overline {1,m} \\ \noa
(\tau^r) \;\succeq \; 0,\quad (\lambda ^\alpha (t))\;\succeq\;0,\;t \in \Omega, \end{array} \right. \leqno \hbox{(VFJ)}
$$
where
$$
\frac{\partial f}{\partial x^0} : = \frac{\partial f}{\partial
x}(t,x^0(t),x_\upsilon ^0 (t)), \qu \frac{\partial f_r }{\partial x_\upsilon ^0}:  = \frac{\partial f_r }{\partial x_\upsilon }(j^1_tx^0)
$$
etc.
\end{Th}

\begin{Def}
The function $x^0(\cdot) \in \mathcal{D}$ is called a
normal efficient solution to (VVP) if, in the conditions (VFJ), there exists $\tau$ with $\tau \succeq 0$, $<e , \tau> = 1$, where $e = (1, \ld 1) \in R^p$.
\end{Def}
\medskip

\subsection{Necessary efficiency for multitime vector \\fractional variational problems (VFP)}

\medskip

In this section we recall some definitions and results that will be needed later in our discussion
about Pareto efficiency conditions for the multitime vector fractional variational problem

$$
\hbox{(VFP)}\; \left\{ \begin{array}{l}
\mbox{Maximize}\; J(x(\cdot)) = \left( \frac{F_1 }{K_1}, \ld, \frac{F_p}{K_p} \right)(x(\cdot)) \\ \noa
\mbox{subject to} \quad g(j^1_tx)\;\preceq\;0, \; h(j^1_tx)\; = \;0, \\ \noa
\quad \quad \quad \quad \quad \quad x(t)\left| {_{\partial \Omega } } \right. = u(t)\;(\mbox{given}),\;\forall t \in \Omega.
\end{array} \right.
$$
The domain of (VFP) is the same set $\mathcal{D}$.

\begin{Def} A feasible solution $x^0(\cdot) \in \mathcal{D}$ is called efficient solution of (VFP)
if there is no $x(\cdot) \in \mathcal{D}$, $x(\cdot) \ne x^0(\cdot)$ such that $J(x(\cdot)) \preceq J(x^0(\cdot))$.
\end{Def}
\medskip

 To present efficiency necessary conditions for (VFP), we need the followings statements.
 Let $x^0(\cdot)$ be an efficient solution to (FVP). Consider the problem
$$
\hbox{(FP)}_r (x^0) \left\{ \begin{array}{l}
 \mathop {\mbox{Minimize}}\limits_x \di\frac{F_r(x(\cdot))}{K_r(x(\cdot))} \\ \noa
\mbox{subject to}\qu x(t)\left| {_{\partial \Omega } =
u(t),\;\forall t \in \Omega } \right. \\ \noa
\quad \quad g(j^1_tx)\preceq 0, \quad
h(j^1_tx)\; = 0 \\ \noa
\quad \quad \di\frac{F_j(x(\cdot))}{K_j(x(\cdot))}\;\leq\di\frac{F_j(x^0(\cdot))}{K_j(x^0(\cdot))}, \; j = \overline {1,p}, \; j \ne r.  \end{array} \right.
$$
Denoting
$$
R_r^0 = \frac{F_r(x^0(\cdot))}{K_r(x^0(\cdot))} = \mathop {\min }\limits_{x(\cdot)} \frac{F_r(x(\cdot))}{K_r(x(\cdot))}, \; r = \overline {1,p},
$$
the problem (FP)$_r (x^0)$ can be written as

$$
\hbox{(FPR)}_r \left\{ {\begin{array}{l}
\mathop {\min }\limits_x \di\frac{F_r(x(\cdot))}{K_r(x(\cdot))} \\ \noa
\mbox{subject to} \qu x(t)\left|_{\partial \Omega} =
u(t),\;\forall t \in \Omega, \right. \\ \noa
g(j^1_tx)\preceq 0, \quad h(j^1_tx) = 0 \\ \noa
F_j (x(\cdot)) - R_j^0 K_j (x(\cdot)) \; \leq 0,\; j = \overline {1,p}, \; j \ne r.
\end{array}} \right.
$$

We add now the problem

$$
\hbox{(SPR)}_r \left\{ {\begin{array}{l}
\mathop {\min }\limits_x (F_r (x(\cdot)) - R_r^0 K_r (x(\cdot)) \\ \noa
\mbox{subject to} \qu (t)\left| {_{\partial \Omega } =
u(t), \; \forall t \in \Omega } \right. \\ \noa
g(j^1_tx)\;\preceq 0,\;h(j^1_tx) = 0 \\ \noa
F_j (x(\cdot)) - R_j^0 K_j (x(\cdot)) \le\;0,\;j = \overline {1,p} ,\;j \ne r. \end{array}} \right.
$$

\begin{Lem}
(Jaganathan [2]). The function $x^0(\cdot) \in \mathcal{D}$ is optimal to (FRP)$_r$ if and only if it is optimal to (SPR)$_r$.
\end{Lem}

\begin{Def} The efficient solution $x^0(\cdot) \in \mathcal{D}$ is called normal efficient solution to (VFP) if
$x^0(\cdot)$ is normal optimal to a least one of scalar problems
(SPR)$_{r\;}$, $r = \overline {1,p}$.\end{Def}

\begin{Th} (Necessary efficiency in (VFP) (Mititelu,
Postolache [8, Theorem 3.2])). Let $x^0(\cdot) \in \mathcal{D}$ be a normal efficient solution to problem (FVP).
Then there exist a vector $\tau = (\tau^r ) \in \textbf{R}^p$ and piecewise smooth functions $\lambda = (\lambda^\alpha (t)) \in \textbf{R}^m$
and $\mu = (\mu^s (t)) \in \textbf{R}^q$ (Lagrange multipliers) that satisfy the conditions
$$
(MFJ) \quad \left\{ \begin{array}{l}
\tau^r \left[ {\di\frac{\partial f_r }{\partial x^0} - R_r^0 \di\frac{\partial k_r }{\partial x^0}} \right] + \lambda^\alpha (t) \di\frac{\partial g_\alpha}{\partial x^0} + \mu^s (t) \di\frac{\partial h_s }{\partial x^0} - \\ \noa
- \di\frac{\partial }{\partial t^\upsilon }\left( {\tau^r \left[
\di{\frac{\partial f_r }{\partial x_\upsilon ^0 } - R_r^0 \di\frac{\partial k_r}{\partial x_\upsilon ^0 }} \right] + \lambda^\alpha (t) \di\frac{\partial g_\alpha }{\partial x_\upsilon^0} + \mu^\beta (t)\frac{\partial h_\beta
}{\partial x_\upsilon^0}} \right) = 0 \\ \noa
\lambda^\alpha (t) g_\alpha (j^1_tx^0) = 0, \qu \hbox{for each} \qu \alpha = \overline {1,m} \\ \noa
\tau \ge 0, \quad <e,\tau> = 1, \; (\lambda^\alpha (t))\;\succeq\;0, \qu t \in \Omega,  \end{array} \right.
$$
where
$$
\frac{\partial f_r }{\partial x^0} = \frac{\partial f_r }{\partial
x}(j^1_tx^0),\;\frac{\partial f_r }{\partial x_\upsilon ^0} = \frac{\partial f_r }{\partial x_\upsilon }(j^1_tx^0)
$$
etc.
\end{Th}

Obviously, we have $R_r^0 = F_r (x^0(\cdot)) / K_r (x^0(\cdot))$, $r = \overline {1,p}$. Taking into account these relations
and denoting $\lambda (t): = K_r (x^0(\cdot))\lambda (t)$, $\mu (t): = K_r (x^0(\cdot))\mu (t)$, Theorem $3$ becomes

\begin{Th} (Necessary efficiency in (VFP) (Mititelu, Postolache [8, Theorem 3.3])). Let $x^0(\cdot) \in \mathcal{D}$ be a normal efficient solution to problem (VFP). Then there exist a vector $\tau = (\tau^r) \in \textbf{R}^p$ and piecewise smooth functions $\lambda = (\lambda^\alpha (t)) \in \textbf{R}^m$ and $\mu = (\mu^\beta (t)) \in  \textbf{R}^q$ (Lagrange multipliers) that satisfy the conditions
$$
(MFJ)_0 \left\{ {\begin{array}{l}
\tau ^r\left[ {K_r (x^0)\di\frac{\partial f_r }{\partial x^0} - F_r
(x^0)\di\frac{\partial k_r }{\partial x^0}} \right] + \lambda ^\alpha (t) \di\frac{\partial g_\alpha }{\partial x^0} + \mu^s (t)\di\frac{\partial h_s }{\partial x^0} - \\ \noa
- \di\frac{\partial }{\partial t^\upsilon }\left( {\tau ^r\left[ {K_r (x^0) \di\frac{\partial f_r }{\partial x_\upsilon ^0 } - F_r (x^0)\di\frac{\partial k_r }{\partial x_\upsilon ^0 }} \right] + \lambda (t)\di\frac{\partial g}{\partial x_\upsilon ^0 } + \mu (t)\di\frac{\partial h}{\partial x_\upsilon^0 }} \right) = 0 \\ \noa  \lambda^\alpha (t) g_\alpha (t,x^0(t),x_\upsilon ^0 (t)) = 0, \qu \hbox{for each} \qu \alpha = \overline {1,m}  \\ \noa
\tau \; \ge \; 0, \quad <{e},\tau> = 1, \; (\lambda^\alpha (t))\;\succeq\;0,\;\;t \in \Omega.  \end{array}} \right.
$$
\end{Th}

\begin{Def}
(Equivalent to Definition $8$). The efficient solution $x^0(\cdot) \in \mathcal{D}$ is called
normal efficient solution to (VFP) if the conditions (MFJ) or (MFJ)$_0 $ exist with $\tau \succeq 0$, $\langle e, \tau \rangle = 1$.
\end{Def}

\section{Sufficient efficiency conditions for \\problems (VVP) and (VFP)}

In this section we shall establish \textit{sufficient efficiency conditions for} (VVP) and (VFP).
In all our theory it is used the essential

{\bf Main Condition}: Suppose that the subset ${\cal S}\subset {\cal F}(\Omega,M)$ is $\eta$-geodesic invex set,
where the $C^1$ vector function $\eta (t)$ is as in Definition $2$, and instead of $\mathcal{F}(\Omega,M)$ we use $\mathcal{S}$.

\begin{Th} (Sufficient efficiency for (VVP)).
Let $x^0(\cdot)\in {\cal S}$, $\tau = (\tau ^r)$, $\lambda = (\lambda ^\alpha )$ and $\mu = (\mu^s)$ be multipliers
satisfying the relations (MFJ) from Theorem 3. For each $x^0(\cdot)\in {\cal S}$, let $x(\cdot)\in {\cal S}$ be an arbitrary geodesic perturbation.
If the following conditions are fulfilled:

a) Each functional $F_r(x(\cdot))$ is $(\rho_r^1 ,b)$-geodesic quasiinvex at $x^0(\cdot)$, with respect to $\eta$ and $d$.

b) The functional $\di\int_\Omega \lambda^\alpha (t)g_\alpha (j^1_tx)\;dv$ is $(\rho_2 ,b)$-geodesic quasiinvex at $x^0(\cdot)$, with respect to $\eta$ and $d$.

c) The functional $\di\int_\Omega \mu^s (t) h_s (j^1_tx)\;dv$ is monotonic $(\rho
_3, b)$-geodesic quasiinvex at $x^0(\cdot)$, with respect to $\eta$ and $d$.

d) One of the functional of a)-b) is strictly $(\rho ,b)$-geodesic quasiinvex at $x^0(\cdot)$.

e) $\tau^r \rho_r^1 + \rho _2 + \rho_3  \geq 0 \quad (\rho_r^1 , \;\rho_2 ,\; \rho_3 \in $ \textbf{R}),

then $x^0(\cdot)$ is an efficient solution to(VVP).
\end{Th}

\textit{Proof}. Let us suppose toward a contradiction, that $x^0(\cdot)$ is not an efficient solution for (VVP). Then, for each $r = \overline {1,p} $,
there exists $x(\cdot) \ne x^0(\cdot)$, a feasible solution to (VVP), such that
\[
\int_\Omega f_r (j^1_tx)\;dv \; \leq \; \int_\Omega f_r (j^1_tx^0)\;dv.
\]

According to hypothesis a), it follows that
\[
b(x,x^0)\int_\Omega [\eta _i (t)\frac{\partial f_r }{\partial
x^i}(j^1_tx^0 ) + \frac{\partial \eta _i }{\partial t^\upsilon}\frac{\partial f_r }{\partial x_\upsilon ^i }(j^1_tx^0 )]\;dv\; \leq \; - \rho_r^1 b(x,x^0)d^2(x,x^0)
\]
(see Definition 2).

Transvecting this inequality by $\tau^r \ge 0$, we obtain
$$
\begin{array}{c}
b(x,x^0)\di\int_\Omega {[\eta _i (t} )\;\tau^r \di\frac{\partial f_r }{\partial x^i}(j^1_tx^0 ) + \di\frac{\partial \eta _i }{\partial t^\upsilon }\tau^r \di\frac{\partial f_r }{\partial x_\upsilon ^i }(j^1_tx^0)]\;dv\;  \\ \noa \leq   - (\tau ^r\rho _r^1 )b(x,x^0)d^2(x,x^0).\end{array} \leqno (4.1)
$$

From the continuity of the functions, we choose $x(t)\in S$ such that
\[
\int_\Omega \lambda^\alpha (t)g_\alpha (t,x,x_\upsilon )\;dv\;\leq   \;\int_\Omega {\lambda ^\alpha } (t)g_\alpha (t,x^0,x_\upsilon^0)\;dv.
\]

Then, taking into account the condition b) and Definition 2, it follows
$$
\begin{array}{c}
b(x,x^0)\di\int_a^b [\eta _i (t) \lambda^\alpha (t) \di\frac{\partial
g_\alpha }{\partial x^i}(j^1_tx^0 ) + \di\frac{\partial \eta _i}{\partial t^\upsilon} \lambda^\alpha (t)\frac{\partial g_\alpha }{\partial x_\upsilon^i}(j^1_tx^0 )] \; dv  \\ \noa \leq   - \rho _2 b(x,x^0)d^2(x,x^0). \end{array} \leqno (4.2)
$$

Taking into account the condition c) and Definition 2, the equality
\[
\int_\Omega \mu^s (t) h_s (j^1_tx)\;dv\; = \; \int_\Omega \mu^s (t) h_s (j^1_tx^0) \;dv
\]
implies
$$
\begin{array}{c}
b(x,x^0) \di\int_\Omega [\eta_i (t)\mu ^s (t) \di\frac{\partial
h_s}{\partial x^i}(j^1_tx^0 )\; + \di\frac{\partial \eta _i }{\partial t^\upsilon }\mu^s (t) \di\frac{\partial h_s }{\partial x_\upsilon^i }(j^1_tx^0 )]\;dv\;  \\ \noa \leq   - \rho _3 b(x,x^0)d^2(x,x^0). \end{array} \leqno (4.3)
$$

We sum side by side the relations (4.1), (4.2) and (4.3) and take into account d). Then we obtain
$$
\begin{array}{l}
b(x,x^0) \di\int_\Omega \eta _i (t) \left[\tau^r \di\frac{\partial f_r }{\partial x^i} + \lambda^\alpha (t) \di\frac{\partial g_\alpha }{\partial x^i} + \mu^s (t) \di\frac{\partial h_s }{\partial x^i} \right]\; (j^1_tx^0 )dv\; + \\ \noa
+ b(x,x^0) \di\int_\Omega \frac{\partial \eta _i }{\partial t^\upsilon } \left[ \tau^r \di\frac{\partial f_r }{\partial x_\upsilon ^i } + \lambda^\alpha (t) \di\frac{\partial g_\alpha }{\partial x_\upsilon ^i } + \mu^s
(t)\frac{\partial h_s }{\partial x_\upsilon ^i}
\right](j^1_tx^0 )\;dv\;  \\ \noa
< - (\tau^r \rho_r^1 + \rho _2 + \rho _3) b(x,x^0)d^2(x,x^0). \end{array} \leqno (4.4)
$$

From (4.4) it follows $b(x,x^0) > 0$ and then
$$
\begin{array}{l}
\di\int_\Omega \eta _i (t) \left[ \tau^r \di\frac{\partial f_r
}{\partial x^i} + \lambda^\alpha (t) \di\frac{\partial g_\alpha }{\partial x^i} + \mu^s (t)\frac{\partial h_s }{\partial x^i} \right]\;(j^1_tx^0 )dv\; + \\ \noa
+ \di\int_\Omega \di\frac{\partial \eta _i }{\partial t^\upsilon } \left[ \tau^r \di\frac{\partial f_r }{\partial x_\upsilon ^i } + \lambda^\alpha (t) \di\frac{\partial g_\alpha }{\partial x_\upsilon ^i } + \mu^s (t) \di \frac{\partial h_s }{\partial x_\upsilon ^i} \right](j^1_tx^0 )\;dv   \\ \noa
< - (\tau ^r\rho _r^1 + \rho _2 + \rho _3 )d^2(x,x^0). \end{array} \leqno (4.5)
$$

We denote $V = [\tau ^rf_r + \lambda ^\alpha (t)g_\alpha + \mu^s (t)h_s ] (j^1_tx^0 )$
and then the relation (4.5) becomes
$$
\int_\Omega {\eta _i } (t)\frac{\partial V}{\partial x^i}dv +
\int_\Omega \frac{\partial \eta _i }{\partial t^\upsilon } \frac{\partial V}{\partial x_\upsilon ^i }dv < - (\tau ^r\rho _r^1 + \rho _2 + \rho _3)d^2(x,x^0), \leqno (4.6)
$$
where we denoted $\frac{\partial V}{\partial x} = \frac{\partial V}{\partial
x}(j^1_tx^0 ),\;\frac{\partial V}{\partial x_\upsilon ^i } = \frac{\partial V}{\partial x_\upsilon ^i }(j^1_tx^0)$.

For an integration by parts in the second integral of (4.5), we have
$$
\frac{\partial }{\partial t^\upsilon }\left( \eta _i (t) \frac{\partial V}{\partial x^i} \right) = \frac{\partial \eta _i}{\partial t^\upsilon} \frac{\partial V}{\partial x^i } + \eta _i (t)\frac{\partial }{\partial t^\upsilon }\left( \frac{\partial V}{\partial x_\upsilon^i}\right), \upsilon = \overline {1,m} .
$$
and
$$
\int_\Omega \frac{\partial \eta _i}{\partial t^\upsilon} \frac{\partial V}{\partial x^i }dv = \int_\Omega \frac{\partial
}{\partial t^\upsilon} \left(\eta _i (t)\frac{\partial V}{\partial x^i} \right)dv - \int_\Omega \eta _i (t)\frac{\partial }{\partial t^\upsilon}\left( \frac{\partial V}{\partial x_\upsilon^i} \right)dv. \leqno (4.7)
$$

Using the flow-divergence formula, we find
$$
\int_\Omega \frac{\partial}{\partial t^\upsilon} \left(\eta _i
(t) \frac{\partial V}{\partial x_\gamma^i} \right) dv = \int_{\partial\Omega} \left( \eta_i (t) \frac{\partial V}{\partial x_\upsilon^i} \right) \; n_\upsilon d\sigma = 0, \leqno (4.8)
$$
where $\bar n = (n_\upsilon )$ is the normal unit vector of the hypersurface $\partial \Omega $ and $\eta (t)\left|_{\partial \Omega} \right. = 0$.
Then, relation (4.7) becomes
$$
\int_\Omega \frac{\partial \eta _i}{\partial t^\upsilon} \frac{\partial V}{\partial x^i}dv = - \int_\Omega \eta_i
(t) \frac{\partial }{\partial t^\upsilon} \left( \frac{\partial V}{\partial x_\upsilon^i} \right)dv \leqno (4.9)
$$
and according to (4.9), the relation (4.6) can be written
$$
\int_\Omega \eta _i (t) \left[\frac{\partial V}{\partial x^i} -
\frac{\partial }{\partial t^\upsilon } \left(\frac{\partial V}{\partial x_\upsilon ^i} \right) \right] dv < - (\tau ^r\rho _r^1 + \rho _2 + \rho_3 )d^2(x,x^0). \leqno (4.10)
$$

Taking into account the first relation of (VVP), the relation (4.10) becomes
$0< - (\tau^r \rho_r^1 + \rho_2 + \rho_3)d^2(x,x^0).$
Having $d(x,x^0) \geq  0$ and the hypothesis e), we obtain the inequality $0 < 0$ that is a false.
Therefore $x^0$ is an efficient solution to (VVP), because $M$ is a complete Riemannian manifold.

In what follows we establish \textit{efficiency sufficient conditions for the problem} (VFP).

\begin{Th}
(Sufficient efficiency for (VFP)). Let $x^0(\cdot)\in {\cal S}$ and $\tau = (\tau^r)$, $\lambda = (\lambda^\alpha)$,
$\mu = (\mu^s )$ be multipliers satisfying the relations(MFJ) from Theorem 3. Suppose fulfilled the following conditions:

a$'$) Each functional $F_r(x(\cdot)) - R_r^0 K_r(x(\cdot))$  is $(\rho_i^1 ,b)$-geodesic quasiinvex at $x^0(\cdot)$,
 with respect to $\eta $ and $d$.

b$'$) The functional $\int_\Omega \lambda^\alpha (t) g_\alpha (j^1_tx)\;dv$ is $(\rho_2 ,b)$-geodesic quasiinvex at $x^0(\cdot)$, with respect to $\eta $ and $d$.

c$'$) The functional $\int_\Omega \mu^s (t) h_s (j^1_tx)\;dv$ is $(\rho_3,b)$-geodesic quasiinvex at $x^0(\cdot)$, with respect to $\eta$ and $d$.

d$'$) One of the functionals from a$'$), b) and c) is strictly $(\rho ,b)$-geodesic quasiinvex at $x^0(\cdot)$, with respect to $\eta $ and $d$ $(\rho = \rho _i^1$, $\rho_2 $ or $\rho_3 $, respectively).

e$'$) $\tau^r \rho_r^1 + \rho _2 + \rho _3 \;\geq\;0$.

Then $x^0(\cdot)$ is an efficient solution to (VFP).
\end{Th}

\textit{Proof}. It is similar to those of Theorem 5, where, for each $r = \overline {1,p} $,
the Lagrangian $f_r (j^1_tx)\;$  is replaced by $f_r (j^1_tx) - R_r^0 k_r (j^1_tx)$.

\begin{Th} (Sufficient efficiency for (VFP)). Let $x^0(\cdot)\in {\cal S}$ and $\tau = (\tau^r)$, $\lambda = (\lambda ^\alpha )$,
$\mu = (\mu^s )$ be multipliers satisfying the relations (MFV)$^0$ from Theorem 4. Suppose 

a$''$) Each functional $r = \overline {1,p}$, $\di\int_\Omega [K_r (x^0) f_r (j^1_tx) -
F_r (x^0)k_r (j^1_tx)]dv$ is $(\rho _r^1 ,b)$-geodesic quasiinvex at $x^0(\cdot)$, with respect to $\eta $ and $d$.

b), c) and e) of Theorem 5.

d$''$) One of the functionals from a$''$), b) and c) is strictly $(\rho ,b)$-geodesic quasiinvex at $x^0(\cdot)$,
with respect to $\eta $ and $d$ ($\rho = \rho _r^1 ,\rho _2 $ or $\rho_3$ respectively).

Then $x^0(\cdot)$ is a geodesic efficient solution to (VFP).
\end{Th}

\textit{Proof}. It is similar to those of Theorem 5, where the hypothesis a) is replaced by hypothesis a$''$) of this theorem.

If, in Theorems 5-7, the functionals from the hypotheses
b) and c) are replaced by the functional $\di\int_\Omega [\lambda^\al (t)g_\al (j^1_tx) + \mu^s(t)h_s (j^1_tx)]dv$, then we have the following results:

\begin{Co}(Sufficient efficiency conditions for (VVP)). Let  $x^0(\cdot)\in {\cal S}$ and $\tau, \lambda$, $\mu$ be
multipliers satisfying the relations (VFJ) from Theorem 3. Suppose 

a) Each functional $r = \overline {1,p}$, $F_r(x(\cdot))$ is $(\rho_r^1 ,b)$-geodesic quasiinvex at $x^0(\cdot)$, with respect to $\eta $ and $d$.

b$'$) The functional $\di\int_\Omega [\lambda^\alpha (t)g_\alpha (j^1_tx) + \mu^s (t) h_s (j^1_tx)] dv$ is $(\rho_2 ,b)$-geodesic quasiinvex at $x^0(\cdot)$, with respect to $\eta $ and $d$.

d$'$) The functionals from a) and b$'$)are strictly $(\rho ,b)$-geodesic quasiinvex at $x^0(\cdot)$, with respect to $\eta $ and $d$ ($\rho = \rho_r^1 $ or $\rho _2 $, respectively).

e$'$) $\tau ^r\rho _r^1 + \rho _2 \; \geq\;0$.

Then $x^0(\cdot)$ is an efficient solution to (VVP).
\end{Co}

\begin{Co} (Sufficient efficiency conditions for (VFP) Let $x^0(\cdot)\in {\cal S}$, and $\tau, \lambda$, $\mu $
be multipliers satisfying the relations (MFJ) from Theorem 3. Suppose satisfied the following conditions:

a$'$) Each functional $r = \overline {1,p}$, $F_r(x(\cdot)) -R_r^0K_r(x(\cdot))$ is $(\rho _r^1 ,b)$-geodesic quasiinvex at $x^0(\cdot)$, with respect to $\eta $ and $d_1$.

b$'$) and e$'$) from Corollary $1$.

d$''$) One of the functionals from a$'$) and b$'$) is strictly $(\rho ,b)$-geodesic quasiinvex at $x^0(\cdot)$, with respect to $\eta $ and $d$.

Then $x^0(\cdot)$ is an efficient solution to (VFP).
\end{Co}

\begin{Co} (Sufficient efficiency conditions for (VFP)). Let $x^0(\cdot)\in {\cal S}$, and $\tau$ ,$\lambda $, $\mu $
be multipliers satisfying the relations (MFV)$_0$ from Theorem 4. Also, we consider
a vector function $\eta $ as in Definitions $3$. Suppose the following conditions are satisfied:

a$''$) Each functional $r = \overline {1,p}$, $\di\int_\Omega [ K_r (x^0)f_r (j^1_tx) - F_r (x^0)k_r (j^1_tx)]\;dv$ is $(\rho _r^1 ,b)$-geodesic quasiinvex at $x^0(\cdot)$, with respect to $\eta $ and $d$.

b$'$) and e$'$) from Corollary $1$.

d$'''$) One of the functionals from a$''$) and b$'$) is strictly $(\rho ,b)$-geodesic quasiinvex at $x^0(\cdot)$, with respect to $\eta $ and $d$, as in Corollary $1$.

Then $x^0(\cdot)$ is an efficient solution to (VFP).
\end{Co}

\begin{Co} (Sufficient optimality conditions for (SVP)). Let $x^0(\cdot)\in {\cal S}$, and $\tau$, $\lambda$, $\mu $
be multipliers satisfying the relations (SFJ) from Theorem 1. If the following conditions 

a) $\di\int_\Omega X(j^1_tx)dv$ is $(\rho ,b)$-geodesic
quasiinvex at $x^0(\cdot)$, with respect to $\eta $ and $d$.

b$'$) From Corollary $1$, c), we have $\tau \rho + \rho _2 \;\ge\;0$.

$\overline d )$ One of the functionals from $\overline a $ and b$'$) is  strictly $(\rho ,b)$-geodesic quasiinvex at $x^0(\cdot)$, with respect to $\eta $ and $d$,

are satisfied, then $x^0(\cdot)$ is an efficient solution to (SVP).
\end{Co}

\bigskip

{\it \c{S}tefan Mititelu

Technical University of Civil Engineering

Balkan Society of Geometry

E-mail: st{\_}mititelu@yahoo.com

\medskip

{\it M\u ad\u alina Constantinescu, Constantin Udri\c{s}te}

University Politehnica of Bucharest

Department of Mathematics and Informatics

313 Splaiul Independen\c tei, 060042 Bucharest

E-mail: udriste@mathem.pub.ro

\end{document}